\documentclass[a4paper,twoside,12pt]{article}
\usepackage[a4paper,margin=3cm,centering,nohead,ignorefoot,footskip=5ex]{geometry}
\usepackage[portuguese]{babel}
\usepackage[T1]{fontenc}
\usepackage[utf8]{inputenc}
\usepackage{amsmath,amssymb,amsthm,url}

\theoremstyle{plain}
\newtheorem*{teoremadacompletude}{Teorema da completude}
\theoremstyle{definition}
\newtheorem*{definicao}{Definição}
\newtheorem*{exemplo}{Exemplo}
\theoremstyle{remark}
\newtheorem*{notacao}{Notação}

\newcommand*{\verdadeiro}{\textnormal V}
\newcommand*{\falso}{\textnormal F}
\newcommand*{\abreviatura}[1]{\vec#1}
\newcommand*{\pontuacao}[1]{\overline{#1}}

\begin{document}

\title{Todas as afirmações verdadeiras são demonstráveis}
\author{Jaime Gaspar\thanks{INRIA Paris-Rocquencourt, $\pi r^2$, Univ Paris Diderot, Sorbonne Paris Cit\'e, F-78153 Le Chesnay, France. Apoiado financeiramente pela Fondation Sciences Math\'ematiques de Paris francesa. \protect\url{mail@jaimegaspar.com}, \protect\url{www.jaimegaspar.com}. Agradeço a Joana Cerveira, Fernando Ferreira, Gilda Ferreira e Rafael Pacheco; no entanto, este texto é da minha inteira responsabilidade.}}
\date{14 de Dezembro de 2012}
\maketitle

\section{Motivação}

Todas as afirmações demonstráveis são verdadeiras. Mas teremos o recíproco: todas as afirmações verdadeiras são demonstráveis? Isto seria o sonho de uma matemática omnisciente capaz de demonstrar todas as afirmações verdadeiras. O teorema da completude dá vida a este sonho:
\begin{equation*}
  \text{todas as \underline{afirmações} \underline{verdadeiras\vphantom{ç}} são \underline{demonstráveis\vphantom{ç}}.} \tag{\textasteriskcentered} \label{teoremadacompletude}
\end{equation*}

Neste texto vamos enunciar e demonstrar o teorema da completude. Como sublinhado em~(\ref{teoremadacompletude}), este teorema fala de afirmações, verdades e demonstrações pelo que, para o podermos enunciar e demonstrar, primeiro precisamos de definir as noções de afirmação, verdade e demonstração.

\section{Afirmação}

Começamos por definir a noção de afirmação. Informalmente, uma afirmação é uma expressão como $\neg(L \vee M) \vee M \vee L$. Recordemos que algumas destas expressões são bem formadas como por exemplo $\neg(L \vee M) \vee M \vee L$, enquanto outras são mal formadas como por exemplo $L \neg\vee M) \vee {\!} \vee ML$.

\begin{definicao}
  Fixemos uma lista de símbolos distintos $L,M,N,\ldots$ aos quais chamamos \emph{letras}. Chamamos \emph{afirmação}, e denotamos por $A,B,C,\ldots$, a uma expressão bem formada construída a partir de letras por meio dos símbolos $\neg$, $\vee$ e dos parênteses.
\end{definicao}

\begin{exemplo}
  A expressão $\neg (L \vee M) \vee M \vee L$ é uma afirmação porque é bem formada e construída a partir de letras $L$ e $M$ por meio de $\neg$, $\vee$ e dos parênteses.
\end{exemplo}

O leitor talvez objete notando que alguns símbolos estão em falta, como por exemplo $\Rightarrow$ e $\exists$. A omissão de $\Rightarrow$ não é relevante: sempre que quisermos falar de $A \Rightarrow B$ podemos em vez disso falar de uma afirmação equivalente que só use $\neg$ e $\vee$ tal como $\neg A \vee B$; por exemplo, no exemplo anterior em vez de falarmos de $L \vee M \Rightarrow M \vee L$, falamos de $\neg(L \vee M) \vee M \vee L$. Já a omissão de $\exists$ é relevante: $\exists x \, A(x)$ não é equivalente a nenhuma afirmação só com $\neg$ e $\vee$. Assim, omitir $\exists$ empobrece de forma relevante as nossas afirmações, mas optámos por o fazer porque simplifica imenso o teorema da completude.

\section{Verdade}

Vamos agora definir a noção de verdade. Denotemos por \verdadeiro{} (respetivamente, \falso) o valor de verdade verdadeiro (respetivamente, falso). Recordemos que podemos calcular o valor de verdade de uma afirmação por meio de uma tabela de verdade. Por exemplo, a tabela de verdade seguinte dá-nos o valor de verdade de $\neg (L \vee M) \vee M \vee L$ em função dos valores de verdade de $L$ e $M$:
\begin{center}
  \begin{tabular}{c|c|c}
    $L$    & $M$    & $\neg (L \vee M) \vee M \vee L$\\\hline
    \verdadeiro  & \verdadeiro  & \verdadeiro\vphantom{$I^{I^{I^I}}$}\\
    \verdadeiro  & \falso & \verdadeiro\\
    \falso & \verdadeiro  & \verdadeiro\\
    \falso & \falso & \verdadeiro
  \end{tabular}
\end{center}

\begin{definicao}
  Dizemos que uma afirmação $A$ é \emph{verdadeira}, e denotamos por $\vDash A$, se e só se todas as entradas da última coluna da tabela de verdade de $A$ são \verdadeiro.
\end{definicao}

\begin{exemplo}
  Temos $\vDash \neg (L \vee M) \vee M \vee L$ porque todas as entradas da última coluna da tabela de verdade de $\neg (L \vee M) \vee M \vee L$ acima são \verdadeiro.
\end{exemplo}

\section{Demonstração}

Finalmente, vamos definir a noção de demonstração. Informalmente, uma demonstração é um argumento em vários passos em que cada passo é um axioma ou resulta de passos anteriores por meio de uma regra. Recordemos que, informalmente: um axioma~$A$ é a asserção ``$A$ é verdadeiro''; uma regra~$\frac{A_1 \ \ \ldots \ \ A_n}{B}$ é a asserção ``se $A_1,\ldots,A_n$ são verdadeiros, então $B$ é verdadeiro'', onde os $A_i$s chamam-se premissas e o $B$ chama-se conclusão.

\begin{notacao}
  Associamos $\smash{\bigvee_{k = 1}^n A_k} = A_1 \vee \cdots \vee A_n$, à direita; por exemplo, $A_1 \vee A_2 \vee A_3 \vee A_4 = A_1 \vee (A_2 \vee (A_3 \vee A_4))$. Denotamos por $\abreviatura A$ uma expressão da forma $\smash{\bigvee_{k = 1}^n A_k}$ (eventualmente $n = 0$, caso em que não há $A_i$s, isto é, $\abreviatura A$ é vazio). Denotamos por $\sigma$ uma permutação de (isto é, uma bijeção de e para) $\{1,\ldots,n\}$.
\end{notacao}

\begin{definicao}
  Consideremos os axiomas e as regras da forma
  \begin{gather*}
    A \vee \neg A \vee \abreviatura B, \tag{A} \label{axioma} \\
    \frac{A_1 \vee \cdots \vee A_n}{A_{\sigma(1)} \vee \cdots \vee A_{\sigma(n)}}, \tag{R$_1$} \label{regra1} \\
    \frac{A \vee (B \vee \abreviatura C)}{(A \vee B) \vee \abreviatura C}, \tag{R$_2$} \label{regra2} \\
    \frac{A \vee \abreviatura B}{\neg\neg A \vee \abreviatura B} \tag{R$_3$} \label{regra3}, \\
    \frac{\neg A \vee \abreviatura C \qquad \neg B \vee \abreviatura C}{\neg(A \vee B) \vee \abreviatura C}. \tag{R$_4$} \label{regra4}
  \end{gather*}
  Dizemos que uma afirmação $A$ é \emph{demonstrável}, e denotamos por $\vdash A$, se e só se existe uma sequência de afirmações $D_1,\ldots,D_n$, chamada \emph{demonstração} de $A$, tal que $A = D_n$, e cada $D_i$ é um dos axiomas, ou é a conclusão de uma das regras sendo a(s) premissa(s) da regra $D_j$(s) com $j < i$.
\end{definicao}

\begin{exemplo}
  Temos $\vdash \neg(L \vee M) \vee M \vee L$ porque o seguinte é uma demonstração de $\neg(L \vee M) \vee M \vee L$:
  \begin{equation*}
    \underbrace{\overbrace{(L \vee M) \vee \neg(L \vee M)}^{D_1}}_{\substack{\text{(\ref{axioma}) com $A = L \vee M$} \\ \text{e $\abreviatura B$ vazio}}}, \ 
    \underbrace{\overbrace{\neg(L \vee M) \vee L \vee M}^{D_2}}_{\substack{\text{conclusão de (\ref{regra1}) com} \\ \text{$n = 2$, $A_1 = L \vee M$,} \\ \text{$A_2 = \neg(L \vee M)$,} \\ \text{$\sigma(1) = 2$ e $\sigma(2) = 1$;} \\ \text{a premissa é $D_1$}}}, \ 
    \underbrace{\overbrace{\neg(L \vee M) \vee M \vee L}^{D_3}}_{\substack{\text{conclusão de (\ref{regra1}) com} \\ \text{$n = 3$, $A_1 = \neg(L \vee M)$,} \\ \text{$A_2 = L$, $A_3 = M$,} \\ \text{$\sigma(1) = 1$, $\sigma(2) = 3$} \\ \text{e $\sigma(3) = 2$;} \\ \text{a premissa é $D_2$}}}.
  \end{equation*}
\end{exemplo}

\section{Teorema da completude}

Depois de termos definido as noções de afirmação, verdade e demonstração, estamos finalmente em condições de enunciar e demonstrar o teorema da completude.

\begin{teoremadacompletude}
  Todas as afirmações verdadeiras são demonstráveis (isto é, para toda a afirmação $A$ tal que $\vDash A$ temos $\vdash A$).
\end{teoremadacompletude}

\begin{proof}
  A cada afirmação $A$ atribuímos uma pontuação $\pontuacao A$ da seguinte forma: cada $\vee$ em $A$ vale $1$ ponto, e cada $\neg$ em $A$ que não esteja imediatamente antes de uma letra vale $1$ ponto. Demonstremos que (1)~$\vDash \smash{\bigvee_{k = 1}^n A_k}$ implica (2)~$\vdash \smash{\bigvee_{k = 1}^n A_k}$ por indução completa em $\smash{\sum_{k = 1}^n \pontuacao{A_k}}$; o teorema é o caso $n = 1$.
  \begin{description}
    \item[Caso base $\smash{\sum_{k = 1}^n \pontuacao{A_k}} = 0$] Temos $\pontuacao{A_i} = 0$ para cada $A_i$, logo cada $A_i$ é uma letra ou a negação de uma letra. Suponhamos (1). Então algum $A_i$ é uma letra $L$ e um algum $A_j$ é $\neg L$, caso contrário dávamos valores de verdade às letras de modo que os $A_k$s fossem falsos contradizendo (1). Seja $\abreviatura A = \smash{\bigvee_{i,j \neq k = 1}^n A_k}$. Temos $\vdash L \vee \neg L \vee \abreviatura A$ por (\ref{axioma}), isto é, $\vdash A_i \vee A_j \vee \abreviatura A$, logo (2) por (\ref{regra1}).

    \item[Passo de indução $\smash{\sum_{k = 1}^n \pontuacao{A_k}} > 0$] Temos $\pontuacao{A_i} > 0$ para algum $A_i$, logo $A_i$ é da forma $B \vee C$ ou $\neg D$, onde $D$ não é uma letra logo é da forma $\neg B$ ou $B \vee C$. Seja $\abreviatura A = \smash{\bigvee_{i \neq k = 1}^n A_k}$. Basta demonstrar que (3)~$\vDash A_i \vee \abreviatura A$ implica (4)~$\vdash A_i \vee \abreviatura A$ porque (1) implica (3) e (4) implica (2) por (\ref{regra1}).
    \begin{description}
      \item[Caso $A_i = B \vee C$] Se (3), isto é, $\vDash (B \vee C) \vee \abreviatura A$, então $\vDash B \vee C \vee \abreviatura A$, logo $\vdash B \vee C \vee \abreviatura A$ por hipótese de indução (que se aplica porque $\pontuacao B + \pontuacao C + \smash{\sum_{i \neq k = 1}^n \pontuacao{A_k}} < \smash{\sum_{k = 1}^n \pontuacao{A_k}}$), portanto $\vdash (B \vee C) \vee \abreviatura A$ por (\ref{regra2}), isto é, (4).

      \item[Caso $A_i = \neg\neg B$] Se (3), isto é, $\vDash \neg\neg B \vee \abreviatura A$, então $\vDash B \vee \abreviatura A$, logo $\vdash B \vee \abreviatura A$ por hipótese de indução (que se aplica porque $\pontuacao B + \smash{\sum_{i \neq k = 1}^n \pontuacao{A_k}} < \smash{\sum_{k = 1}^n \pontuacao{A_k}}$), portanto $\vdash \neg\neg B \vee \abreviatura A$ por (\ref{regra3}), isto é, (4).

      \item[Caso $A_i = \neg (B \vee C)$] Se (3), isto é, $\vDash \neg(B \vee C) \vee \abreviatura A$, então $\vDash \neg B \vee \abreviatura A$ e $\vDash \neg C \vee \abreviatura A$, logo $\vdash \neg B \vee \abreviatura A$ e $\vdash \neg C \vee \abreviatura A$ por hipótese de indução (que se aplica porque $\pontuacao{\neg B} + \smash{\sum_{i \neq k = 1}^n \pontuacao{A_k}},\, \pontuacao{\neg C} + \smash{\sum_{i \neq k = 1}^n \pontuacao{A_k}} < \smash{\sum_{k = 1}^n \pontuacao{A_k}}$), portanto $\vdash \neg(B \vee C) \vee \abreviatura A$ por (\ref{regra4}), isto é, (4). \cite{Srivastava2001}\qedhere
    \end{description}
  \end{description}
\end{proof}

\section{Sugestões de leitura}

Se o leitor estiver interessado em saber mais, o passo seguinte é ler sobre:
\begin{description}
  \item[Teorema da completude de G\"odel] Estende o teorema aqui apresentado de modo a abranger afirmações com $\exists$. \cite{DetlovsPodnieks2000,Kennedy2007,Srivastava2001}
  \item[Teoremas da incompletude de G\"odel] Demonstram que existem afirmações verdadeiras (num sentido mais geral do que o usado neste texto) que são indemonstráveis. \cite{Kennedy2007,Podnieks1997,Srivastava2007}
\end{description}

\bibliography{Referencias}{}

\begin{thebibliography}{1}

\bibitem{DetlovsPodnieks2000}
Vilnis~Detlovs e~Karlis~Podnieks.
\newblock Introduction to mathematical logic, 2000.
\newblock \url{http://www.ltn.lv/~podnieks/mlog/ml.htm}.

\bibitem{Kennedy2007}
Juliette Kennedy.
\newblock {K}urt {G}{\"{o}}del.
\newblock In Edward~N. Zalta, editor, {\em Stanford Encyclopedia of
  Philosophy}. Fevereiro 2007.
\newblock \url{http://plato.stanford.edu/entries/goedel/}.

\bibitem{Podnieks1997}
Karlis Podnieks.
\newblock What is mathematics: {G}{\"{o}}del's theorem and around, 1997.
\newblock \url{http://www.ltn.lv/~podnieks/gt.html}.

\bibitem{Srivastava2001}
S.~M. Srivastava.
\newblock The completeness theorem of {G}{\"{o}}del.
\newblock {\em Resonance}, 6(7, 8), Julho, Agosto 2001.
\newblock \url{http://www.ias.ac.in/resonance/}.

\bibitem{Srivastava2007}
S.~M. Srivastava.
\newblock {G}{\"{o}}del's proof.
\newblock {\em Resonance}, 12(2, 3, 5), Fevereiro, Março, Maio 2007.
\newblock \url{http://www.ias.ac.in/resonance/}.

\end{thebibliography}
\bibliographystyle{plain}

\end{document}